\title{\vspace{-0.7cm} Two remarks on the Burr-Erd\H{o}s conjecture}
\author{
Jacob Fox\thanks{Department of Mathematics, Princeton, Princeton,
NJ. Email: {\tt jacobfox@math.princeton.edu}. Research supported by
an NSF Graduate Research Fellowship and a Princeton Centennial
Fellowship.} \and Benny Sudakov\thanks{Department of Mathematics,
UCLA,  Los Angeles, CA 90095. Email: {\tt bsudakov@math.ucla.edu}.
Research supported in part by NSF CAREER award DMS-0546523, NSF
grant DMS-0355497, and by a USA-Israeli BSF grant.}}

\documentclass[11pt]{article}
\usepackage{amsfonts,amssymb,amsmath,latexsym}

\newenvironment{proof}
      {\medskip\noindent{\bf Proof.}\hspace{1mm}}
      {\hfill$\Box$\medskip}

\def\qed{\ifvmode\mbox{ }\else\unskip\fi\hskip 1em plus 10fill$\Box$}

\newtheorem{theorem}{Theorem}[section]

\newtheorem{lemma}[theorem]{Lemma}

\newtheorem{corollary}[theorem]{Corollary}

\newtheorem{definition}[theorem]{Definition}

\begin{document}
\date{}

\maketitle

\begin{abstract}

The Ramsey number $r(H)$ of a graph $H$ is the minimum positive
integer $N$ such that every two-coloring of the edges of the
complete graph $K_N$ on $N$ vertices contains a monochromatic copy
of $H$. A graph $H$ is $d$-degenerate if every subgraph of $H$ has
minimum degree at most $d$. Burr and Erd\H{o}s in 1975 conjectured
that for each positive integer $d$ there is a constant $c_d$ such
that $r(H) \leq c_dn$ for every $d$-degenerate graph $H$ on $n$
vertices. We show that for such graphs $r(H) \leq 2^{c_d\sqrt{\log
n}}n$, improving on an earlier bound of Kostochka and Sudakov. We
also study Ramsey numbers of random graphs, showing that for $d$
fixed, almost surely the random graph $G(n,d/n)$ has Ramsey number
linear in $n$. For random bipartite graphs, our proof gives nearly
tight bounds.
\end{abstract}

\section{Introduction}

For a graph $H$, the {\it Ramsey number} $r(H)$ is the least
positive integer $N$ such that every two-coloring of the edges of
complete graph $K_N$ on $N$ vertices contains a monochromatic copy
of $H$. Ramsey's theorem states that $r(H)$ exists for every graph
$H$. A classical result of Erd\H{o}s and Szekeres, which is a
quantitative version of Ramsey's theorem, implies that $r(K_n) \leq
2^{2n}$ for every positive integer $n$. Erd\H{o}s showed using
probabilistic arguments that $r(K_n) > 2^{n/2}$ for $n
> 2$. Over the last sixty years, there has been several
improvements on these bounds (see, e.g., \cite{Co}). However,
despite efforts by various researchers, the constant factors in the
above exponents remain the same.

Determining or estimating Ramsey numbers is one of the central
problem in combinatorics, see the book {\it Ramsey theory}
\cite{GrRoSp} for details. Besides the complete graph, the next most
classical topic in this area concerns the Ramsey numbers of sparse
graphs, i.e., graphs with certain upper bound constraints on the
degrees of the vertices. The study of these Ramsey numbers was
initiated by Burr and Erd\H{o}s in 1975, and this topic has since
placed a central role in graph Ramsey theory.

A graph is {\it $d$-degenerate} if every subgraph has a vertex of
degree at most $d$. In 1975, Burr and Erd\H{o}s \cite{BuEr}
conjectured that, for each positive integer $d$, there is a constant
$c(d)$ such that every $d$-degenerate graph $H$ with $n$ vertices
satisfies $r(H) \leq c(d)n$. An important special case of this conjecture for
bounded degree graphs was proved by Chv\'atal, R\"odl, Szemer\'edi,
and Trotter \cite{ChRoSzTr}.

Another notion of sparseness was introduced by Chen and Schelp
\cite{ChSc}. A graph is {\it $p$-arrangeable} if there is an
ordering $v_1,\ldots,v_n$ of the vertices such that for any vertex
$v_i$, its neighbors to the right of $v_i$ have together at most $p$
neighbors to the left of $v_i$ (including $v_i$). This is an
intermediate notion of sparseness not as strict as bounded degree
though not as general as bounded degeneracy. Extending the result of
\cite{ChRoSzTr}, Chen and Schelp proved that there is a constant
$c(p)$ such that every $p$-arrangeable graph $H$ on $n$ vertices has
Ramsey number at most $c(p)n$. This gives linear Ramsey numbers for
planar graphs and more generally for graphs that can be drawn on a
bounded genus surface. This result was later extended by R\"odl and
Thomas \cite{RoTh}, who showed that graphs with no $K_p$-subdivision
are $p^8$-arrangeable.

Here we introduce a notion of sparseness that is closely related
to arrangeability. The main reason for introducing this notion
is that it turns out to be more useful for bounding Ramsey numbers. A graph $H$ is {\it
$(d,\Delta)$-degenerate} if there exists an ordering
$v_1,\ldots,v_n$ of its vertices such that for each $v_i$,
\begin{enumerate}
\item there are at most $d$ vertices $v_j$ adjacent to $v_i$ with
$j<i$, and
\item there are at most $\Delta$ subsets $S \subset
\{v_1,\ldots,v_i\}$ such that $S=N(v_j) \cap \{v_1,\ldots,v_i\}$ for some neighbor $v_j$ of $v_i$ with $j>i$,
where the {\it neighborhood} $N(v_j)$ is the set of vertices that are adjacent to $v_j$.
\end{enumerate}

From the definition, every $(d,\Delta)$-degenerate graph is
$d$-degenerate, and every graph with maximum degree $\Delta$ is
$(\Delta,\Delta)$-degenerate. More interesting but also very simple
to show is that every $(d,\Delta)$-degenerate graph is
$\left(\Delta(d-1)+1\right)$-arrangeable, and every $p$-arrangeable
graph is $(p,2^{p-1})$-degenerate (see Lemmas \ref{arrangedeg} and
\ref{arrangedeg1}).

While the conjecture of Burr and Erd\H{o}s is still open, there has been considerable progress on this problem recently.
Kostochka and R\"odl \cite{KoRo1} were the first to prove a
polynomial upper bound on the Ramsey numbers of $d$-degenerate
graphs. They showed that $r(H) \leq c_dn^2$ for every $d$-degenerate
graph $H$ with $n$ vertices. A nearly linear bound of the form $r(H)
\leq 2^{c_d(\log n)^{2d/(2d+1)}}n$ was obtained by Kostochka and Sudakov
\cite{KoSu}. In \cite{FoSu1}, the authors proved that $r(H) \leq 2^{c_d\sqrt{\log n}}n$ for every bipartite $d$-degenerate graph $H$ with
$n$ vertices. Here we show how to use the techniques developed in \cite{FoSu} to generalize this result to all $d$-degenerate graphs.

\begin{theorem}\label{main}
For each positive integer $d$ there is a constant $c_d$
such that every $(d,\Delta)$-degenerate graph $H$ with order $n$ satisfies
$r(H) \leq 2^{c_d\sqrt{\log \Delta}} \, n.$ In particular, $r(H) \leq 2^{c_d \sqrt{\log n}} \, n$ for every $d$-degenerate graph $H$
on $n$ vertices.
\end{theorem}

\noindent This result follows from Theorem \ref{ramseydegenerategeneral}, which gives a more general bound
on the Ramsey number which also incorporates the chromatic number of the graph $H$.

We next discuss Ramsey numbers and arrangeablity of sparse random graphs. The {\it random graph} $G(n,p)$ is the probability space of labeled graphs on $n$ vertices, where every edge appears
independently with probability $p=p(n)$. We say that the random graph possesses a graph property $\mathcal{P}$
 {\it almost surely}, or a.s.\ for brevity, if the probability that $G(n,p)$ has property $\mathcal{P}$ tends to $1$ as
 $n$ tends to infinity. It is well known and easy to show that if $p=d/n$ with $d >0$ fixed, then, a.s.\ $G(n,d/n)$ will
have maximum degree $\Theta(\log n/\log \log n)$. Moreover, as shown by
Ajtai, Koml\'os, and Szemer\'edi \cite{AKS}, for $d>1$ fixed, a.s.\ the random graph $G(n,d/n)$ contains
a subdivision of $K_p$ with $p$ almost as large as its maximum degree. Therefore one cannot use known results to give a
linear bound on the Ramsey number of $G(n,d/n)$. Here we obtain such a bound by proving that the random graph $G(n,d/n)$ a.s.\
has bounded arrangeability.

\begin{theorem}\label{randomarrange}
There are constants $c_1>c_2>0$ such that for $d \geq 1$ fixed,
a.s.\ the random graph $G(n,d/n)$ is $c_1d^2$-arrangeable but not $c_2d^2$-arrangeable.
\end{theorem}

\noindent Theorem \ref{randomarrange} is closely related to a question of Chen and Schelp, who asked to estimate the proportion of
$d$-degenerate graphs which have bounded arrangeability. The following is an immediate corollary of this theorem and the result of Chen and Schelp.

\begin{corollary}\label{corrandom} For each $d \geq 1$, there is a constant $c_d$ such that a.s.\ the Ramsey number of
$G(n,d/n)$ is at most $c_dn$.
\end{corollary}

One would naturally like to obtain good bounds on the Ramsey number
of $G(n,d/n)$. To accomplish this, we prove a stronger version of
Theorem \ref{randomarrange}, showing that for $d \geq 300$ a.s.\
$G(n,d/n)$ is $(16d,16d)$-degenerate. We also modify a proof of
Graham, R\"odl, Ruci\'nski, to prove that there is a constant $c$ such
that for every $(d,\Delta)$-degenerate $H$ with $n$ vertices and
chromatic number $q$, $r(H) \leq \left(d^d\Delta\right)^{c\log q}n$.
In the other direction, it follows from a result of Graham, R\"odl,
and Ruci\'nski \cite{GrRoRu} that a.s.\ $G(n,d/n)$ has Ramsey number
at least $2^{cd}n$ for some absolute positive constant $c$. From
these results we get the following quantitative version of Corollary
\ref{corrandom}.

\begin{theorem}\label{corrandom1}
There are positive constants $c_1,c_2$ such that for $d \geq 2$ and
$n$ a.s.\
$$2^{c_1d}\, n \leq r(G(n,d/n)) \leq 2^{c_2d\log^2 d} \, n.$$
\end{theorem}

In the case of random bipartite graphs, we can obtain nearly tight
bounds. In another paper, Graham, R\"odl, and Ruci\'nski
\cite{GrRoRu1} adapt their proof of a lower bound for Ramsey numbers
to work also for random bipartite graphs. The random bipartite graph
$G(n,n,p)$ is the probability space of labeled bipartite graphs with
$n$ vertices in each class, where each of the $n^2$ edges appears
independently with probability $p$.

\begin{theorem}\label{corrandom2}
There are positive constants $c_1,c_2$ such that for each $d \geq 1$
and $n$ a.s.\
$$2^{c_1d} \, n \leq r(G(n,n,d/n)) \leq 2^{c_2d} \, n.$$
\end{theorem}

The rest of this paper is organized as follows. In the next section
we present a proof of Theorem \ref{ramseydegenerategeneral} which
implies Theorem \ref{main} on the Ramsey number for
$(d,\Delta)$-degenerate graphs. In Section \ref{next}, we prove
another bound on Ramsey numbers for $(d,\Delta)$-degenerate graphs
which is sometimes better than Theorem
\ref{ramseydegenerategeneral}. In Section \ref{randomnext}, we prove
results on the random graphs $G(n,d/n)$ and $G(n,n,d/n)$. The last
section of this paper contains some concluding remarks. Throughout
the paper, we systematically omit floor and ceiling signs whenever
they are not crucial for the sake of clarity of presentation. We
also do not make any serious attempt to optimize absolute constants
in our statements and proofs. All logarithms in this paper are base
$2$.

\section{Proof of Theorem \ref{main}}

The main result of this section is the following general bound on
Ramsey numbers.

\begin{theorem}\label{ramseydegenerategeneral}
There is a constant $c$ such that for $0<\delta \leq 1$, every
$(d,\Delta)$-degenerate graph $H$ with chromatic number $q$ and
order $n$ satisfies
$$r(H)<2^{cq3^qd/\delta}\Delta^{c\delta}n.$$
\end{theorem}

Note that a greedy coloring shows that every $d$-degenerate graph
has chromatic number at most $d+1$. Theorem \ref{main} follows from
the above theorem by letting $\delta=1/\sqrt{\log \Delta}$.

The first result that we need is a lemma from \cite{FoSu1} whose
proof uses a probabilistic argument known as {\it dependent random
choice}. Early versions of this technique were developed in the
papers \cite{Go2,KoRo,Su}. Later, variants were discovered and
applied to various Ramsey-type problems (see, e.g.,
\cite{KoSu,AlKrSu,Su1,FoSu1}, and their references). We include the
proof here for the sake of completeness. Given a vertex subset $T$
of a graph $G$, the {\it common neighborhood $N(T)$ of $T$} is the
set of all vertices of $G$ that are adjacent to $T$, i.e., to {\it
every} vertex in $T$.

\begin{lemma}\label{lem:dependent}
If $\epsilon>0$ and $G=(V_1,V_2;E)$ is a bipartite graph with
$|V_1|=|V_2|=N$ and at least $\epsilon N^2$ edges, then for all
positive integers $t,x$, there is a subset $A \subset V_2$ with $|A|
\geq \epsilon^{2t}N/2$ such that for all but at most
$2\epsilon^{-2t}\left(\frac{x}{N}\right)^{2t}{N \choose t}$ $t$-sets
$S$ in  $A$, we have $|N(S)| \geq x$.
\end{lemma}
\begin{proof}
Let $T$ be a subset of $2t$ random vertices of $V_1$, chosen
uniformly with repetitions. Set $A=N(T)$, and let $X$ denote the
cardinality of $A \subset V_2$. By linearity of expectation and by
convexity of $f(z)=z^{2t}$,
$$\mathbb{E}[X]=\sum_{v \in V_2}\left(\frac{|N(v)|}{N}\right)^{2t}=
N^{-2t}\sum_{v \in V_2}|N(v)|^{2t} \geq N^{1-2t}\left(\frac{\sum_{v
\in V_1}|N(v)|}{N}\right)^{2t} \geq \epsilon^{2t}N.$$

Let $Y$ denote the random variable counting the number of $t$-sets
in $A$ with fewer than $x$ common neighbors. For a given $t$-set
$S$, the probability that $S$ is a subset of $A$ is
$\left(\frac{|N(S)|}{N}\right)^{2t}$. Therefore, we have
$$\mathbb{E}[Y] \leq {N \choose t}\left(\frac{x-1}{N}\right)^{2t}.$$
Using linearity of expectation, we have
$$\mathbb{E}\left[X-\frac{\mathbb{E}[X]}{2\mathbb{E}[Y]}\, Y-\mathbb{E}[X]/2\right] \geq 0.$$
Therefore, there is a choice of $T$ for which this expression is
non negative. Then
$$|A|=X \geq \frac{1}{2}\mathbb{E}[X] \geq \frac{1}{2}\epsilon^{2t}N$$
and $$Y\leq 2X\mathbb{E}[Y]/\mathbb{E}[X] \leq
2N\mathbb{E}[Y]/\mathbb{E}[X]<
2\epsilon^{-2t}\left(\frac{x}{N}\right)^{2t}{N \choose t},$$
completing the proof.
\end{proof}

We use this lemma to deduce the following:

\begin{lemma}\label{firststep}
For every $2$-edge-coloring of $K_N$ and integers $t \geq q \geq 2$
there is a color and nested subsets of vertices $A_1 \subset \ldots
\subset A_q$ with $|A_1| \geq 2^{-4tq}N$ such that the following
holds. For each $i<q$, all but at most  $2^{4t^2q}y^{2t}N^{-t}$
subsets of $A_{i}$ of size $t$ have at least $y$ common neighbors in
$A_{i+1}$ in this color.
\end{lemma}
\begin{proof}
For $j \in \{0,1\}$, let $G_j$ denote graph of color $j$. Let
$B_1=V(K_N)$. We will pick subsets $B_1 \supset B_2 \supset \ldots
\supset B_{2q-2}$ such that for each $i \in [2q-3]$, we have
$|B_{i+1}| \geq |B_i|/2^{2t+2}$ and there is a color $c(i) \in
\{0,1\}$ such that there are less than $\left(8y^2/|B_i|\right)^t$
$t$-sets $S \subset B_{i+1}$ which have less than $y$ common
neighbors in $B_i$ in graph $G_{c(i)}$.

Having already picked $B_i$, we now show how to pick $c(i)$ and
$B_{i+1}$ . Arbitrarily partition $B_i$ into two subsets $B_{i,1}$
and $B_{i,2}$ of equal size. Let $c(i)$ denote the densest of the
two colors between $B_{i,1}$ and $B_{i,2}$. By Lemma
\ref{lem:dependent} with $\epsilon=1/2$, there is a subset $B_{i+1}
\subset B_{i,2} \subset B_i$ with $|B_{i+1}| \geq 2^{-2t-1}|B_{i,2}|
=2^{-2t-2}|B_i|$ such that for all but at most
$$2 \cdot
2^{2t}\left(\frac{y}{|B_{i,2}|}\right)^{2t}{|B_{i,2}| \choose t}\leq
2^{3t}y^{2t}|B_i|^{-t}=\left(8y^2/|B_i|\right)^t$$ $t$-sets $S
\subset B_{i+1}$, $S$ has at least $y$ common neighbors in $B_i$ in
graph $G_{c(i)}$.

We have completed the part of the proof where we constructed the
nested subsets $B_1 \supset \ldots \supset B_{2q-2}$ and the colors
$c(1),\ldots,c(2q-3)$. Notice that $|B_{2q-2}| \geq
2^{-(2t+2)(2q-3)}N \geq 2^{-4tq+6}N$. So for all but at most
$$\left(8y^2/|B_i|\right)^t \leq
\left(2^{4tq-3}y^2/N\right)^t \leq 2^{4t^2q}y^{2t}N^{-t}$$ $t$-sets
$S \subset B_{i+1}$, $S$ has at least $y$ common neighbors in $B_i$
by color $c(i)$. Since the sets $B_1 \supset \ldots \supset
B_{2q-2}$ are nested, this also implies that that for all but at
most $2^{4t^2q}y^{2t}N^{-t}$ $t$-sets $S \subset B_{i+1}$, $S$ has
at least $y$ common neighbors in $B_j$ in graph $G_{c(i)}$ for each
$j \leq i$.

By the pigeonhole principle, one of the two colors is represented at
least $q-1$ times in the sequence $c(1),\ldots,c(2q-3)$. We suppose
without loss of generality that $0$ is this popular color. Let
$A_q=B_1$, $i_j$ denote the $j^{\textrm{th}}$ smallest positive
integer such that $c(i_j)=0$, and $A_{q-j}=B_{i_j+1}$ for $1 \leq j
\leq q-1$. By the above discussion, it follows that $A_1 \subset
\ldots A_{q}$, $|A_1| \geq |B_{2q-2}| \geq 2^{-4t^2q}N$, and, for
each positive integer $i<q$, all but at most $2^{4t^2q}y^{2t}N^{-t}$
subsets of $A_{i}$ of size $t$ are adjacent to at least $y$ vertices
in $A_{i+1}$ in graph $G_0$. This completes the proof.
\end{proof}

The previous lemma shows that in every $2$-edge-coloring of the
complete graph $K_N$ there is a monochromatic subgraph $G$ and large
nested vertex subsets $A_1 \subset \ldots \subset A_q$ such that
almost every $t$-set in $A_i$ has large common neighborhood in
$A_{i+1}$ in graph $G$. The next lemma is the most technical part of
the proof. It says that if a graph $G$ has such vertex subsets $A_1
\subset \ldots \subset A_q$, then for $1 \leq i \leq q$ there are
large subsets $V_i \subset A_i$ such that almost every $d$-set in
$\bigcup_{\ell \not = i}V_{\ell}$ has large common neighborhood in
$V_i$. These vertex subsets $V_1,\ldots,V_q$ will be used to show
that $G$ contains all $(d,\Delta)$-degenerate graphs on $n$ vertices
with chromatic number at most $q$.

\begin{lemma}\label{mainstep} Let $d,q,\Delta \geq 2$ be integers and $0<\delta \leq 1$.
Let $t=(3^{q}-1)d/\delta+d$, $y=2^{-5qt}\Delta^{-\delta}N$, and
$x=y^4N^{-3}$. Suppose $x \geq 2t$. Let $G=(V,E)$ be a graph with
nested vertex subsets $A_1 \subset \ldots \subset A_q$ with
$|A_1|\geq 2^{-4tq}N$ such that for each $i$, all but at most
$2^{4t^2q}y^{2t}N^{-t}$ subsets of $A_{i}$ of size $t$ have
 at least $y$ common neighbors in $A_{i+1}$. Then there are vertex subsets $V_i \subset A_i$ for $1 \leq i \leq q$ such that
 $|V_i| \geq x$ and the number of $d$-sets in $\bigcup_{\ell \not = i}V_{\ell}$
 with fewer than $x$ common neighbors in $V_i$ is less than $(2\Delta)^{-d}{x \choose d}$.
\end{lemma}

\medskip\noindent{\bf Proof.}\hspace{1mm}
We will first pick some constants. Let $r_0=t$ and for $1 \leq j
\leq q$, let $t_j=2 \cdot 3^{q-j}d/\delta$ and
$r_{j}=r_{j-1}-t_j=(3^{q-j}-1)d/\delta+d$. In particular, we have
$r_q=d$ and $t_j \geq 2r_j$.

Let $$b_{i,i}=2q\left(2x/y\right)^{t_i}{N \choose
r_i}~~~~\textrm{and}~~~~
b_{i,j}=2q\left(\frac{r_{j-1}}{y}\right)^{t_j}b_{i,j-1}~~\textrm{for}~i<j.$$

%Note that
%\begin{eqnarray*}
%B_j & = & B_1\prod_{i=2}^j B_i/B_{i-1}=
%2q \left(2x/y\right)^{t_1}{N \choose r_1}\prod_{i=2}^j 2q\left(\frac{r_{j-1}}{y}\right)^{t_j}
%= (2q)^{j}(2t^{-10})^{t_1}\left(\frac{y}{N}\right)^{t_1}{N \choose r_1}\prod_{i=2}^j \left(\frac{r_{j-1}}{y}\right)^{t_j} \\
%& \geq &
%.
%\end{eqnarray*}

%We have
%\begin{eqnarray*}
%B_q &=&(2q)^{q-1} \left(\prod_{j=2}^{q}
%\left(2\frac{r_{j-1}}{y}\right)^{t_j}\right)B_1 \leq (2q)^{q-1}(2t)^t
%y^{d-r_1}B_1 \leq t^{2t}y^{d-r_1-t_1}x^{t_1}{N \choose r_1} \\ &\leq& t^{2t}
%\left(\frac{x}{y}\right)^{t_1-d}\left(\frac{N}{y}\right)^{r_1}\frac{x^d}{r_1!}
%=
%t^{2t-10t_1+10d}\left(\frac{N}{y}\right)^{r_1-3t_1+3d}\frac{x^d}{r_1!}<\left(\frac{N}{y}\right)^{-2t_1}\frac{x^d}{r_1!}<(2\Delta)^{-d}{x
%\choose d}.
%\end{eqnarray*}
Let
$$c_0=2^{-t^2q}y^t~~~~\textrm{and}~~~~c_j=2q\left(\frac{r_{j-1}}{y}\right)^{t_j}c_{j-1}.$$

By the hypothesis of the lemma, we have nested subsets $A_1 \subset
\ldots \subset A_q$ with $|A_1|\geq 2^{-4tq}N$ such that for each
$i$, all but at most $2^{4t^2q}y^{2t}N^{-t} \leq 2^{-t^2q}y^t=c_0$
subsets of $A_{i}$ of size $t$ are adjacent to
 at least $y$ vertices in $A_{i+1}$. Let $A_{i,0}=A_{i}$ for each $i$. We will prove by induction on $j$
that there are subsets $A_{1,j},\ldots,A_{q,j}$ for $1 \leq j \leq
q$ that satisfy the following properties.

\begin{enumerate}
\item For $1 \leq i,j \leq q$ , $A_{i,j} \subset A_{i,j-1}$.
\item For $0 \leq j<\ell<i \leq q$, $A_{\ell,j} \subset A_{i,j}$.
\item $|A_{i,j}| \geq y$ for all $i>j$
and $|A_{i,j}| \geq 2x-t_i$ for $i \leq j$.
\item For each $i \leq j$, the number of $r_j$-sets in $\bigcup_{\ell \not = i} A_{\ell,j}$ that have
less than $2x-t_i$ common neighbors in $A_{i,j}$ is at most
$b_{i,j}$.
\item For $j<i<q$, the number of $r_j$-sets in $A_{i,j}$ with less than $y$ common neighbors in $A_{i+1,j}$ is at most
$c_j$.
\end{enumerate}

It is easy to see that the desired properties hold for $j=0$. Assume
we have already found the subsets $A_{i,j-1}$ for $1 \leq i \leq q$.
We now show how to pick the subsets $A_{i,j}$. Pick a subset $S_j$
of $A_{j,j-1}$ of size $t_j$ uniformly at random. We will let
$A_{i,j}=A_{i,j-1} \cap N(S_j)$ for $i \not =j$ and
$A_{j,j}=A_{j,j-1}\setminus S_j$.

Let $X_j$ be the random variable that counts the number of
$r_j$-sets in $\bigcup_{\ell \not = j} A_{\ell,j}$ with at most
$2x-t_j$ common neighbors in $A_{j,j}$. Equivalently, $X_j$ is the
number of $r_j$-sets in $\bigcup_{\ell \not = j} A_{\ell,j}$ with at
most $2x$ common neighbors in $A_{j,j-1}$. The number of $r_j$-sets
in $\bigcup_{\ell \not = j} A_{\ell,j-1}$ is at most ${N \choose
r_j}$, and the probability that a given $r_j$-set $R \subset
\bigcup_{\ell \not = j} A_{\ell,j-1}$ with at most $2x$ common
neighbors in $A_{j,j-1}$ is also contained in $\bigcup_{\ell \not =
j} A_{\ell,j}$ is at most ${2x \choose t_j}\big/ {|A_{j,j-1}|
\choose t_j}$. By linearity of expectation, we have
$$\mathbb{E}[X_j] \leq \frac{{2x \choose t_j}}{{|A_{j,j-1}| \choose
t_j}}{N \choose r_j} \leq \left(\frac{2x}{y}\right)^{t_j}{N \choose
r_j} = b_{j,j}/2q.$$

For $i < j$, let $Y_{i,j}$ be the random variable that counts the
number of $r_{j-1}$-sets containing $S_j$ in $\bigcup_{\ell \not =
i} A_{\ell,j-1}$ with less than $2x-t_i$ common neighbors in
$A_{i,j-1}$. Since
 the number of $r_{j-1}$-sets in $\bigcup_{\ell \not = i}
A_{\ell,j-1}$ that have less than $2x-t_i$ common neighbors in
$A_{i,j-1}$ is at most $b_{i,j-1}$, then
$$\mathbb{E}[Y_{i,j}] \leq \frac{{r_{j-1} \choose t_j}}{{|A_{j,j-1}| \choose t_j}}b_{i,j-1} \leq \left(\frac{r_{j-1}}{y}\right)^{t_j}b_{i,j-1} =b_{i,j}/2q.$$

Note that $A_{\ell,j}$ is disjoint from $S_j$ for each $\ell$. So if
$T$ is a subset of $\bigcup_{\ell \not = i} A_{\ell,j}$ with
cardinality $r_j$, then $T$ is disjoint from $S_j$ and so $T \cup
S_j$ is a subset of $\bigcup_{\ell \not = i} A_{\ell,j-1}$ with
cardinality $r_j+t_j=r_{j-1}$ satisfying
$$|N(T \cup
S_j) \cap A_{i,j-1}|=|N(T) \cap N(S_j) \cap A_{i,j-1}|=|N(T) \cap
A_{i,j}|.$$ Hence, $Y_{i,j}$ is also an upper bound on the number of
$r_j$-sets in $\bigcup_{\ell \not = i} A_{\ell,j}$ with less than
$2x-t_i$ common neighbors in $A_{i,j}$.

For $j<i <q$, let $Z_{i,j}$ be the random variable that counts the
number of $r_{j-1}$-sets in $A_{i,j-1}$ that contain $S_j$ and have
less than $y$ common neighbors in $A_{i+1,j-1}$. Since $c_{j-1}$ is
an upper bound on the number of $r_{j-1}$-sets in $A_{i,j-1}$ that
have less than $y$ common neighbors in $A_{i+1,j-1}$ and
 $S_j \subset A_{j,j-1} \subset A_{i,j-1}$,
 then
$$\mathbb{E}[Z_{i,j}] \leq \frac{{r_{j-1} \choose t_j}}{{|A_{j,j-1}| \choose t_j}}c_{j-1} \leq \left(\frac{r_{j-1}}{y}\right)^{t_j}c_{j-1} =c_j/2q.$$
Since $A_{i,j}$ is disjoint from $S_j$, if $T \subset A_{i,j}$ has
cardinality $r_j$, then $T \cup S_j \subset A_{i,j-1}$ has
cardinality $r_j+t_j=r_{j-1}$ and satisfies
$$|N(T \cup S_j) \cap A_{i+1,j-1}|=|N(T) \cap N(S_j) \cap A_{i+1,j-1}|=|N(T) \cap
A_{i+1,j}|.$$ Hence, $Z_{i,j}$ is an upper bound on the number of
subsets of $A_{i,j}$ of size $r_{j}$ with less than $y$ common
neighbors in $A_{i+1,j}$.

For $i<j$, let $F_{i,j}$ be the event that every $r_{j-1}$-set in
$A_{j,j-1}$ containing $S_j$ has less than $2x-t_i$ common neighbors
in $A_{i,j-1}$. The number of $r_{j-1}$-sets in $A_{j,j-1}$ with
less than $2x-t_i$ common neighbors in $A_{i,j-1}$ is at most
$b_{i,j-1}$. The number of subsets of $A_{j,j-1}$ of size $r_{j-1}$
containing a fixed subset of size $t_j$ is ${|A_{j,j-1}|-t_j \choose
r_{j-1}-t_j}$ and there are ${r_{j-1} \choose t_j}$ subsets of size
$t_j$ in an $r_{j-1}$-set. Hence, the number of $t_j$-sets in
$A_{j,j-1}$ for which every $r_{j-1}$-set in $A_{j,j-1}$ containing
it has less than $2x-t_i$ common neighbors in $A_{i,j-1}$ is at most
$$\frac{{r_{j-1} \choose t_j}}{{|A_{j,j-1}|-t_j \choose
r_{j-1}-t_j}}b_{i,j-1}.$$ Since there are a total of ${|A_{j,j-1}|
\choose t_j}$ possible $t_j$-sets in $A_{j,j-1}$ that can be picked
for $S_j$ and $|A_{j,j-1}| \geq y$, then the probability of event
$F_{i,j}$ is at most

\begin{eqnarray*}\frac{{r_{j-1} \choose t_j}}{{|A_{j,j-1}|-t_j \choose
r_{j-1}-t_j}}b_{i,j-1}{|A_{j,j-1}| \choose t_j}^{-1} & \leq &
\frac{r_{j-1}^{t_j}}{t_j!} \cdot
\frac{(r_{j-1}-t_j)!}{(y/2)^{r_{j-1}-t_j}} \cdot
\frac{t_j!}{(y/2)^{t_j}} \cdot b_{i,j-1} =
2^{r_{j-1}}r_{j-1}^{t_j}r_j!y^{-r_{j-1}} b_{i,j-1}
\\ & \leq &
t^{t_i}y^{-r_{j-1}}b_{i,j-1} = t^{t_i}y^{-r_{j-1}}b_{i,i}\prod_{i
\leq \ell < j-1}\frac{b_{i,\ell+1}}{b_{i,\ell}} \\ & = &
t^{t_i}y^{-r_{j-1}}2q\left(2x/y\right)^{t_i}{N \choose r_i}\prod_{i
\leq \ell < j-1} 2q\left(\frac{r_{\ell}}{y}\right)^{t_{\ell+1}}
\\ & \leq & t^{t_i}y^{-r_{j-1}}(2q)^{q}\left(2x/y\right)^{t_i}N^{r_i}\left(\frac{r_i}{y}\right)^{t_{i+1}+\cdots+t_{j-1}}
\\ & = & (2t)^{t_i}(2q)^{q}\left(x/y\right)^{t_i}\left(N/y\right)^{r_i}r_i^{r_i-r_{j-1}}
< t^{3t_i}(2q)^{q}\left(x/y\right)^{t_i}\left(N/y\right)^{r_i} \\
& = & t^{3t_i}(2q)^{q}\left(N/y\right)^{r_i-3t_i}\leq
t^{3t_i}(2q)^{q}\left(2^{5tq}\right)^{-5t_i/2} <\frac{1}{2q},
\end{eqnarray*}
where we used $t=r_0 \geq r_1 \ldots \geq r_q$, $r_{j-1}=t_j+r_j$
and $t_j \geq 2r_j$ for $1 \leq j \leq q$.

If there is an $r_{j-1}$-set $T$ in $A_{j,j-1}$ containing $S_j$
with at least $2x-t_i$ common neighbors in $A_{i,j-1}$, then $S_j$
has at least $2x-t_i$ common neighbors in $A_{i,j-1}$ and
$|A_{i,j}|=|A_{i,j-1} \cap N(S_j)| \geq 2x-t_i$. Therefore, if
$|A_{i,j}|<2x-t_i$, then event $F_{i,j}$ occurs.

Let $G_{j}$ be the event that every $r_{j-1}$-set in $A_{j,j-1}$
containing $S_j$ has less than $y$ common neighbors in
$A_{j+1,j-1}$. The number of $r_{j-1}$-sets in $A_{j,j-1}$ that have
less than $y$ common neighbors in $A_{j+1,j-1}$ is at most
$c_{j-1}$. The number of subsets of $A_{j,j-1}$ of size $r_{j-1}$
containing a fixed set of size $t_j$ is ${|A_{j,j-1}|-t_j \choose
r_{j-1}-t_j}$ and there are ${r_{j-1} \choose t_j}$ subsets of size
$t_j$ in an $r_{j-1}$-set. Hence, the number of $t_j$-sets in
$A_{j,j-1}$ for which every $r_{j-1}$ set in $A_{j,j-1}$ containing
it has less than $y$ common neighbors in $A_{j+1,j-1}$ is at most
$$\frac{{r_{j-1} \choose t_j}}{{|A_{j,j-1}|-t_j \choose
r_{j-1}-t_j}}c_{j-1}.$$ Since there are a total of ${|A_{j,j-1}|
\choose t_j}$ possible $t_j$-sets in $A_{j,j-1}$ that can be picked
for $S_j$ and $|A_{j,j-1}| \geq y$, then the probability of event
$G_{i,j}$ is at most

\begin{eqnarray*}
\frac{{r_{j-1} \choose t_j}}{{|A_{j,j-1}|-t_j \choose
r_{j-1}-t_j}}c_{j-1}{|A_{j,j-1}| \choose t_j}^{-1} & \leq &
\frac{r_{j-1}^{t_j}}{t_j!} \cdot
\frac{(r_{j-1}-t_j)!}{(y/2)^{r_{j-1}-t_j}} \cdot
\frac{t_j!}{(y/2)^{t_j}} \cdot c_{j-1} =
2^{r_{j-1}}r_{j-1}^{t_j}r_j!y^{-r_{j-1}} c_{j-1}
\\ & \leq & t^{2t}y^{-r_{j-1}}c_{j-1} = t^{2t}y^{-r_{j-1}}c_{0}\prod_{\ell=1}^{j-1}
\frac{c_{\ell}}{c_{\ell-1}} \\ & = & t^{2t}y^{-r_{j-1}} 2^{-t^2q}y^t
\prod_{\ell=1}^{j-1} 2q\left(\frac{r_{\ell-1}}{y}\right)^{t_{\ell}}
=t^{2t}2^{-t^2q}(2q)^{j-1} \prod_{\ell=1}^{j-1}
\left(r_{\ell-1}\right)^{t_{\ell}} \\ & < &
t^{2t}2^{-t^2q}(2q)^{j-1}t^{t} < t^{4t}2^{-t^2q} < \frac{1}{2q},
\end{eqnarray*}
where we used $t=r_0 \geq r_1 \ldots \geq r_q$ and
$r_{j-1}=t_j+r_j$.

If there is an $r_{j-1}$-set $T$ in $A_{j,j-1}$ containing $S_j$
with at least $y$ common neighbors in $A_{j+1,j-1}$, then $S_j$ has
at least $y$ common neighbors in $A_{j+1,j-1}$ and
$|A_{j+1,j}|=|A_{j+1,j-1} \cap N(S_j)| \geq y$. Therefore, if
$|A_{j+1,j}|<y$, then event $G_j$ occurs.

Note that each of the discrete random variables $X_j,Y_{i,j},
Z_{i,j}$ are non negative. Markov's inequality for non negative random
variables says that if $X$ is a non negative discrete random variable
and $c \geq 1$, then the probability that $X >  c\mathbb{E}[X]$ is
less than $\frac{1}{c}$. So each of the following five types of
events have probability less than $\frac{1}{2q}$ of occurring:
\begin{enumerate}
\item $F_{i,j}$ with $i<j$,
\item $G_j$,
\item $X_j > 2q\mathbb{E}[X_j]$,
\item $Y_{i,j} > 2q\mathbb{E}[Y_{i,j}]$ with $i<j$,
\item $Z_{i,j} > 2q\mathbb{E}[Z_{i,j}]$ with $j<i$.
\end{enumerate}
Since there are a total of $j-1+1+1+j-1+q-j=q+j\leq 2q$ events of
the above five types, there is a positive probability that none of
these events occur. Hence, there is a choice of $S_j$ for which none
of the events $F_{i,j}$ occur, the event $G_j$ does not occur,
$$X_j \leq 2q\mathbb{E}[X_j] \leq b_{j,j},~\hspace{0.2cm}
Y_{i,j} \leq 2q\mathbb{E}[Y_{i,j}]  \leq
b_{j,j}~\textrm{for}~i<j,~\hspace{0.2cm}~\textrm{and}~
\hspace{0.2cm}~ Z_{i,j} \leq 2q\mathbb{E}[Z_{i,j}]  \leq
c_j~\textrm{for}~j<i.$$

Recall that $A_{i,j} =A_{i,j-1} \cap N(S_j)$ if $i \not = j$ and
$A_{j,j}=A_{j,j-1} \setminus S_j$. Hence $A_{i,j} \subset A_{i,j-1}$
for $1 \leq i,j\leq q$, which is the first of the five desired
properties. By the induction hypothesis, $A_{\ell,j-1} \subset
A_{i,j-1}$ for $j-1<\ell<i$, and so for $j<\ell<i$,
$A_{\ell,j}=A_{\ell,j-1} \cap N(S_j) \subset A_{i,j-1} \cap
N(S_j)=A_{i,j}$, which is the second of the desired properties
follows.

Note that $|A_{i,j}| \geq 2x-t_i$ for $i<j$ since $F_{i,j}$ does not
occur, and $|A_{j,j}| =|A_{j,j-1}|-t_j \geq y-t_j \geq 2x-t_j$. For
$i>j$, since $A_{j+1,j} \subset \ldots \subset A_{q,j}$ and event
$G_j$ does not occur, then $|A_{i,j}| \geq |A_{j+1,j}|\geq y$. This
demonstrates the third of the five desired properties. The upper
bounds on $X_j$ and the $Y_{i,j}$ shows the fourth desired property
and the upper bound on the $Z_{i,j}$ shows the fifth desired
property. Hence, by induction on $j$, the $A_{i,j}$ have the desired
properties.

We let $V_i=A_{i,q}$ for $1 \leq i \leq q$. For each $i$, we have
$|V_i| \geq 2x-t_i \geq x$ and all but less than $b_{i,q}$ $d$-sets
in $\bigcup_{\ell \not = i} V_{\ell}$ have at least $x$ common
neighbors in $V_{i}$. To complete the proof, it suffices to show
that $b_{i,q} < (2\Delta)^{-d}{x \choose d}$. Using $t=r_0 \geq r_1
\ldots \geq r_q=d$, $r_{i-1}=t_i+r_i$ and $t_i \geq 2r_i$ for $1
\leq i \leq q$, we have
\begin{eqnarray*}b_{i,q} & = & b_{i,i}\prod_{i
\leq \ell \leq q-1}\frac{b_{i,\ell+1}}{b_{i,\ell}} =
2q\left(2x/y\right)^{t_i}{N \choose r_i}\prod_{i \leq \ell \leq q-1}
2q\left(\frac{r_{\ell}}{y}\right)^{t_{\ell+1}} \leq
(2q)^{q}\left(2x/y\right)^{t_i}N^{r_i}
\left(\frac{r_{i}}{y}\right)^{t_{i+1}+\cdots+t_q} \\ & = &
(2q)^{q}\left(2x\right)^{t_i}y^{r_q-r_{i-1}}N^{r_i}r_{i}^{t_{i+1}+\cdots+t_q}
= (2q)^q2^{t_i}r_i^{r_i-d}x^{t_i}y^{d-r_{i-1}}N^{r_i}\leq
2^{qtt_i}x^{t_i}y^{d-r_{i-1}}N^{r_i}
\\ & = &  2^{qtt_i}x^d\left(y^4N^{-3}\right)^{t_i-d}y^{d-r_{i-1}}N^{r_i}
 =  2^{qtt_i}x^{d}(y/N)^{3t_i-3d-r_i} =
2^{qtt_i}x^{d}(2^{-5qt}\Delta^{-\delta} )^{3t_i-3d-r_i}
\\ & \leq & 2^{qtt_i}x^{d}2^{-5qtt_i}\Delta^{-2d} = 2^{-4qtt_i}\Delta^{-2d}x^d <  (2\Delta)^{-d}{x \choose d}.
\hspace{6.5cm} \Box
\end{eqnarray*}

We now present the proof of Theorem \ref{ramseydegenerategeneral}.
Given a $2$-coloring of the edges of $K_N$ with $N \geq 2^{25q
 3^q d/\delta}\Delta^{4\delta}n$ we must show that it
contains a monochromatic copy of every $(d,\Delta)$-degenerate graph
$H$ with $n$ vertices and chromatic number $q$. Using Lemmas
\ref{firststep} and \ref{mainstep} (with
$y=2^{-5tq}\Delta^{-\delta}N$ and $t \leq 3^qd/\delta$), we can find
vertex subsets $V_1,\ldots,V_q$, each of cardinality at least
$$x=(y/N)^4N=(2^{-5qt}\Delta^{-\delta})^4N \geq 2^{5qt}n> 4n,$$ and a
monochromatic subgraph $G$ of the $2$-edge-coloring such that for
each $i$ the number of $d$-sets $S \subset \bigcup_{\ell \not =
i}V_{\ell}$ with fewer than $x$ common neighbors in $V_i$ is less
than $(2\Delta)^{-d}{x \choose d}$. Then the following embedding
lemma shows that $G$ contains a copy of $H$ and completes the proof.

\begin{lemma}\label{degeneracyhelpful}
Let $H$ be a $(d,\Delta)$-degenerate graph with $n$ vertices and
chromatic number $q$. Let $G$ be a graph with vertex subsets
$V_1,\ldots, V_q$ such that for each $i$, $|V_i| \geq x \geq 4n$ and
the number of $d$-sets $S \subset \bigcup_{\ell \not = i}V_{\ell}$
with fewer than $x$ common neighbors in $V_i$
 is less than $(2\Delta)^{-d}{x \choose d}$. Then $G$
contains a copy of $H$.
\end{lemma}
\begin{proof}
Call a $d$-set $S \subset \bigcup_{\ell \not = i} V_{\ell}$ {\it
good with respect to $i$} if $|N(S)\cap V_i| \geq x$, otherwise it
is {\it bad with respect to $i$}. Also, a subset $S \subset
\bigcup_{\ell \not = i} V_{\ell}$ with $|S| < d$ is {\it good with
respect to $i$} if it is contained in less than $(2\Delta)^{|S|-d}{x
\choose d-|S|}$ $d$-sets in $\bigcup_{\ell \not = i} V_{\ell}$ which
are bad with respect to $i$. A vertex $v \in V_{k}$ with $k \not =
i$ is {\it bad with respect to $i$ and a subset $S \subset
\bigcup_{\ell \not = i} V_{\ell}$} with $|S|<d$ if $S$ is good with
respect to $i$ but $S \cup \{v\}$ is not. Note that, for any subset
$S \subset \bigcup_{\ell \not = i} V_{\ell}$ with $|S|<d$ that is
good with respect to $i$, there are at most $\frac{x}{2\Delta}$
vertices that are bad with respect to $S$ and $i$. Indeed, if not,
then there would be more than
$$\frac{x/(2\Delta)}{d-|S|}(2\Delta)^{|S|+1-d}{x \choose d-|S|-1}\geq (2\Delta)^{|S|-d}{x\choose d-|S|}$$
subsets of $\bigcup_{\ell \not = i} V_{\ell}$ of size $d$ containing
$S$ that are bad with respect to $i$, which would contradict $S$
being good with respect to $i$.

Since $H$ is $(d,\Delta)$-degenerate, it has an ordering of its
vertices $\{v_1,\ldots,v_n\}$ such that each vertex $v_k$ has at
most $d$ neighbors $v_{\ell}$ with $\ell<k$ and there are at most
$\Delta$ subsets $S \subset \{v_1,\ldots,v_{k}\}$ such that
$S=N(v_j) \cap \{v_1,\ldots,v_k\}$ for some neighbor $v_j$ of $v_k$
with $j>k$. Since $H$ has chromatic number $q$, there is a partition
$U_1 \cup \ldots \cup U_q$ of the vertex set of $H$ such that each
$U_i$ is an independent set. For $1 \leq j \leq n$, let $r(j)$
denote the index $r$ of the independent set $U_r$ containing vertex
$v_j$. Let $N^-(v_k)$ be all the neighbors $v_{\ell}$ of $v_{k}$
with $\ell<k$. Let $L_h=\{v_1,\ldots,v_h\}$. An {\it embedding} of a
graph $H$ in a graph $G$ is an injective mapping $f:V(H) \rightarrow
V(G)$ such that $(f(v_j),f(v_k))$ is an edge of $G$ if $(v_j,v_k)$
is an edge of $H$. In other words, an embedding $f$ demonstrates
that $H$ is a subgraph of $G$. We will use induction on $h$ to find
an embedding $f$ of $H$ in $G$ such that for $1 \leq i \leq q$,
$f(U_i) \subset V_i$ and for every vertex $v_k$ and every $h \in
[n]$, the set $f(N^-(v_k)\cap L_h)$ is good with respect to $r(k)$.

By our definition, the empty set is good with respect to each $i$,
$1 \leq i \leq q$. We will embed the vertices of $H$ by increasing
order of their indices. Suppose we are embedding $v_h$. Then, by the
induction hypothesis, for each vertex $v_k$, the set $f(N^-(v_k)\cap
L_{h-1})$ is good with respect to $r(k)$. Since the set
$f(N^-(v_h)\cap L_{h-1})=f(N^-(v_h))$ is good with respect to
$r(h)$, it has at least $x$ common neighbors in $V_{r(h)}$. Also,
there are at most $\Delta$ subsets $S \subset L_h$ for which there
is a neighbor $v_j$ of $v_h$ with $j>h$ such that $L_h \cap
N(v_j)=S$. So there are at most $\Delta$ sets $f(N^-(v_j) \cap
L_{h-1})$ where $v_j$ is a neighbor of $v_h$ and $j>h$. Note that
$r(j) \not = r(h)$ for $v_j$ a neighbor of $v_h$ with $j>h$. By the
induction hypothesis each such set $f(N^-(v_j) \cap L_{h-1})$ is
good with respect to $r(j)$, so there are at most $\Delta
\frac{x}{2\Delta} = x/2$ vertices in $V_{r(h)}$ which are bad with
respect to at least one of the pairs $f(N^-(v_j) \cap L_{h-1})$ and
$r(j)$. This implies that there at least $x-x/2-(h-1)>x/4$ vertices
in $V_{r(h)}$ in the common neighborhood of $f(N^-(v_h))$ which are
not occupied yet and are good with respect to all the above pairs
$f(N^-(v_j) \cap L_{h-1})$ and $r(j)$. Any of these vertices can be
chosen as $f(v_h)$. When the induction is complete, $f(v_h)$ is
adjacent to $f(N^-(v_h))$ for every vertex $v_h$ of $H$. Hence, the
mapping $f$ provides an embedding of $H$ as a subgraph of $G$.
\end{proof}
\section{Another bound for Ramsey numbers}\label{next}

The following theorem is a generalization of a bound by Graham,
R\"odl, and Ruci\'nski \cite{GrRoRu} on Ramsey numbers for graphs of
bounded maximum degree. The proof is a minor variation on their
proof.

\begin{theorem}\label{GRRtheorem}
The Ramsey number of every $(d,\Delta)$-degenerate $H$ with $n$
vertices and chromatic number $q$ satisfies $r(H) \leq
\left(2^{7d+8}d^{3d+2}\Delta\right)^{\log q}n$.
\end{theorem}

\noindent We will need the following lemma. The {\it edge density} between a pair of vertex subsets $W_1,W_2$ of a graph $G$ is
the fraction of pairs $(w_1,w_2) \in W_1 \times W_2$ that are edges of $G$.

\begin{lemma}\label{GRR}
Let $\epsilon>0$ and $H$ be a $(d,\Delta)$-degenerate graph with $n$
vertices and chromatic number $q$. Let $G$ be a graph on $N \geq
4\epsilon^{-d}qn$ vertices. If every pair of disjoint subsets
$W_1,W_2 \subset V(G)$ each with cardinality at least
$\frac{1}{2}\epsilon^{d}\Delta^{-1}q^{-2}N$
 has edge density at least $\epsilon$ between them, then $G$ contains $H$ as a subgraph.
\end{lemma}
\begin{proof}
Let $v_1,\ldots,v_n$ be a $(d,\Delta)$-degenerate ordering for $H$.
Let $L_j=\{v_1,\ldots,v_j\}$. Let $V(H)=U_1 \cup \ldots \cup U_q$ be
a partition of the vertex set of $H$ into $q$ color classes which
are independent sets. For $i
>j$, let $N(i,j)=N(v_i) \cap L_j$ denote the set of neighbors $v_h$ of $v_i$ with
$h \leq j$ and $d_{i,j}=|N(i,j)|$. Arbitrarily partition $V(G)=V_1
\cup \ldots \cup V_q$ into $q$ subsets of size $N/q$. We will find
an embedding $f:V(H) \rightarrow V(G)$ of $H$ such that if $v_i \in
U_k$, then $f(v_i) \in V_k$. For $1 \leq i \leq n$ and $v_i \in
U_k$, let $T_{i,0}=V_k$. We will embed the vertices $v_1,\ldots,v_n$
of $H$ one by one in increasing order. We will prove by induction on
$j$ that at the end of step $j$, we will have vertices
$f(v_1),\ldots,f(v_j)$ and sets $T_{i,j}$ for $i>j$ such that
$$|T_{i,j}| \geq \epsilon^{d_{i,j}}|T_{i,0}| = \epsilon^{d_{i,j}}N/q
\geq \epsilon^d N/q$$ and the following holds. For $h,\ell \leq j$,
$(f(v_h),f(v_{\ell}))$ is an edge of $G$ if $(v_h,v_{\ell})$ is an
edge of $H$, and for $i>j$ and $v_i \in U_k$, $T_{i,j}$ is the
subset of $V_k$ consisting of those vertices adjacent to $f(v_p)$
for every vertex $v_p \in N(i,j)$. In particular, we have that if
$i,i'>j$ are such that $v_i,v_{i'}$ lie in the same independent set
$U_k$ and $N(i,j)=N(i',j)$, then $T_{i,j}=T_{i',j}$. Note that any
vertex in $T_{i,j} \setminus \{f(v_1),\ldots,f(v_j)\}$ can be used
to embed $v_i$.

The base case $j=0$ for the induction clearly holds. Our induction
hypothesis is that we have vertices $f(v_1),\ldots,f(v_{j-1})$ and
sets $T_{i,j-1}$ for $i>j-1$ with $|T_{i,j-1}| \geq
\epsilon^{d_{i,j-1}}N/q$ such that if $(v_h,v_{\ell})$ is an edge of
$H$ and $h,\ell<j$, then $\left(f(v_{h}),f(v_{\ell})\right)$ is an
edge of $G$, and if $(v_h,v_{\ell})$ is an edge of $H$ with $h<j
\leq \ell$, then $f(v_h)$ is adjacent to every element of
$T_{\ell,j-1}$. It is sufficient to find a vertex $w \in T_{j,j-1}
\setminus \{f(v_1),\ldots,f(v_{j-1})\}$ such that for each $v_i$
adjacent to $v_j$ with $i>j$, the number of elements of $T_{i,j-1}$
adjacent to $w$ is at least $\epsilon|T_{i,j-1}|$. Indeed, if we
find such a vertex $w$, we let $f(v_j)=w$ and for $i>j$, we let
$T_{i,j}=N(w) \cap T_{i,j-1}$ if $v_i$ is adjacent to $v_j$ in $H$
and otherwise $T_{i,j}=T_{i,j-1}$, which completes step $j$. For
$v_i$ adjacent to $v_j$ with $i>j$, let $X_{i,j}$ denote the set of
vertices in $T_{j,j-1}$ with less than $\epsilon|T_{i,j-1}|$
neighbors in $T_{i,j-1}$. If there is a $X_{i,j}$ with cardinality
at least $\frac{1}{2q\Delta}|T_{j,j-1}|$, then letting $W_1=X_{i,j}$
and $W_2=T_{i,j-1}$, the edge density between $W_1$ and $W_2$ is
less than $\epsilon$ and $W_1,W_2$ each have cardinality at least
$\frac{1}{2}\epsilon^{d}\Delta^{-1}q^{-2}N$, contradicting the
assumption of the lemma. So each of the sets $X_{i,j}$ have
cardinality less than $\frac{1}{2q\Delta}|T_{j,j-1}|$.

Since $H$ is $(d,\Delta)$-degenerate, there are at most $\Delta$
vertex subsets $S \subset L_j$ with $v_j \in S$ for which there is a
vertex $v_i$ with $i>j$ and $N(i,j)=S$. As we already mentioned, if
$i,i'>j-1$ are such that $v_i,v_{i'}$ lie in the same independent
set $U_k$ and $N(i,j-1)=N(i',j-1)$, then $T_{i,j-1}=T_{i',j-1}$. If
furthermore $v_i,v_{i'}$ are neighbors of $v_j$, then
$X_{i,j}=X_{i',j}$. Since there are $q$ sets $U_k$ and at most
$\Delta$ sets $S \subset L_j$ with $v_j \in S$ for which there is a
vertex $v_i$ with $i>j$ and $N(i,j)=S$, then there are at most
$q\Delta$ distinct sets of the form $X_{i,j}$. Hence, at least
$$|T_{j,j-1}|-q\Delta\frac{1}{2q\Delta}|T_{j,j-1}|-(j-1)>\frac{1}{2}|T_{j,j-1}|-n
\geq \frac{1}{2}\epsilon^d q^{-1}N-n
\geq n$$ vertices $w \in T_{j,j-1}$ can be chosen for $f(v_j)$,
which by induction on $j$ completes the proof.
\end{proof}

We now mention some useful terminology from \cite{FoSu} that we need
before proving Theorem \ref{GRRtheorem}. For a graph $G=(V,E)$ and
disjoint subsets $W_1,\ldots,W_t \subset V$, the {\it density}
$d_{G}(W_1,\ldots,W_t)$ between the $t \geq 2$ vertex subsets
$W_1,\ldots,W_t$ is defined by
$$d_G(W_1,\ldots,W_t)=\frac{\sum_{ i < j}
e(W_i,W_j)}{\sum_{i < j } |W_i||W_j|}.$$ If
$|W_1|=\ldots=|W_t|$, then
$$d_G(W_1,\ldots,W_t)={t \choose 2}^{-1}\sum_{ i < j }
d_G(W_i,W_j).$$

\begin{definition} \label{d31}
For $\alpha,\rho,\epsilon \in [0,1]$ and positive integer $t$, a graph $G=(V,E)$ is
{\bf $(\alpha,\rho,\epsilon,t)$-sparse} if for all subsets $U \subset V$ with $|U| \geq \alpha |V|$,
there are disjoint subsets $W_{1},\ldots,W_{t} \subset U$ with $|W_{1}|=\ldots =|W_{t}|=\lceil \rho |U|\rceil$
and $d_{G}(W_1,\ldots,W_{t_i})\leq \epsilon$.
\end{definition}

By averaging, if $\alpha'\geq \alpha$, $\rho' \leq \rho$,
$\epsilon'\geq \epsilon$, $t' \leq t$, and $G$ is
$(\alpha,\rho,\epsilon,t)$-sparse, then $G$ is also
$(\alpha',\rho',\epsilon',t')$-sparse. To prove Theorem
\ref{GRRtheorem}, we use the following simple lemma from
\cite{FoSu}.

\begin{lemma} \cite{FoSu} \label{FoSulemma} If $G$ is $(\alpha,\rho,\epsilon/4,2)$-sparse, then $G$
is also
$\left((\frac{2}{\rho})^{h-1}\alpha,2^{1-h}\rho^{h},\epsilon,2^h\right)$-sparse for each positive integer $h$.
\end{lemma}

For this paragraph, let $\epsilon=\frac{1}{32q^2d}$,
$x=4\epsilon^{-d}\Delta q^2$, and $y=(2x)^{\log q}$. Note that Lemma
\ref{GRR} demonstrates that if a graph $G$ on $N \geq xn$ vertices
does not contain a $(d,\Delta)$-degenerate graph $H$ with order $n$
and chromatic number $q$, then $G$ is
$(xn/N,x^{-1},\epsilon,2)$-sparse. By Lemma \ref{FoSulemma} with
$\alpha=xn/N$, $h=\log q$, $\rho=1/x$, this implies that if a graph
$G$ on $N \geq yn$ vertices does not contain a
$(d,\Delta)$-degenerate graph $H$ with order $n$ and chromatic
number $q$, then $G$ is $\left(yn/N,y^{-1},
4\epsilon,q\right)$-sparse. Hence, as long as $N \geq yn$, then
there are vertex subsets $W_1,\ldots,W_q$ of $G$ with the same size
which is at least $N/y$ such that $d_G(W_1,\ldots,W_q)$ is at most
$4\epsilon$. Consider a red-blue edge-coloring of $K_N$ with
$$N= \left(2^{7d+8}d^{3d+2}\Delta\right)^{\log q}n \geq
8\left(8q^2(32dq^2)^{d}\Delta\right)^{\log
q}n=8\left(8\epsilon^{-d}\Delta q^2\right)^{\log q}n=8yn,$$ where we
use the fact that the chromatic number $q$ of a $d$-degenerate graph
satisfies $q \leq d+1 \leq 2d$. If the red graph does not contain
$H$, then there are disjoint subsets $W_1,\ldots,W_q$ of $V(K_N)$
each with the same cardinality which is at least $y^{-1}N \geq 8n$
such that $d_{G(R)}(W_1,\ldots,W_q)$ is at most $4\epsilon$, where
$G(R)$ denotes the graph of color red. Hence, the total number of
red edges whose vertices are in different $W_i$s is at most
$4\epsilon q^2|W_1|^2$. For each $W_i$, delete those $|W_i|/2$
vertices of $W_i$ which have the largest number of neighbors in
$\bigcup_{j \not = i} W_j$ in the red graph, and let $Y_i$ be the
remaining vertices of $W_i$. Notice that $|Y_i|=|W_i|/2 \geq 4n$ and
every vertex of $Y_i$ is in at most $8\epsilon q^2
|W_i|=\frac{|Y_i|}{2d}$ red edges with vertices in $\bigcup_{j \not
= i} W_j$ since otherwise the number of edges between $W_i \setminus
Y_i$ and $\bigcup_{j \not =i} W_j$ is more than $|W_i \setminus
Y_i|8 \epsilon  q^2 |W_i| = 4\epsilon q^2|W_1|^2$, contradicting the
fact that the number $\sum_{i<j} e(W_i,W_j)$ of edges with vertices
in different subsets is at most $4\epsilon q^2 |W_1|^2$. Therefore,
applying the following lemma to the blue graph, there is a
monochromatic blue copy of $H$, completing the proof of
\ref{GRRtheorem}.

\begin{lemma}
Suppose $H$ is a $d$-degenerate graph with $n$ vertices and
chromatic number $q$. If $G$ is a graph with disjoint vertex subsets
$Y_1,\ldots,Y_q$ with $|Y_1|=\ldots=|Y_q| \geq 4n$ such that each
vertex in $Y_i$ is adjacent to all but at most $\frac{|Y_i|}{2d}$
vertices of $\bigcup_{j \not = i}Y_j$, then $H$ is a subgraph of
$G$.
\end{lemma}
\begin{proof}
Let $v_1,\ldots,v_n$ be an ordering of the vertices of $H$ such that
for each vertex $v_i$, there are at most $d$ neighbors $v_j$ of
$v_i$ with $j<i$. Let $V(H)=U_1 \cup \ldots \cup U_q$ be a partition
of the vertex set of $H$ into independent sets. We will find an
embedding $f:V(H) \rightarrow V(G)$ of $H$ such that if $v_i \in
U_k$, then $f(v_i) \in Y_k$. We will embed the vertices
$v_1,\ldots,v_n$ of $H$ one by one in increasing order. Suppose we
have already embedded $v_1,\ldots,v_{i-1}$ and we try to embed $v_i
\in U_k$. Consider a vertex $v_j$ adjacent to $v_i$ with $j<i$. We
have $v_j \not \in U_k$ since $U_k$ is an independent set.
Therefore, $f(v_j) \not \in Y_k$ and is adjacent to all but at most
$\frac{|Y_k|}{2d}$ vertices in $Y_k$. Since there are at most $d$
such vertices $v_j$ adjacent to $v_i$ with $j<i$, then there are at
least $|Y_k|-d\frac{|Y_k|}{2d}-(i-1)=\frac{|Y_k|}{2}-(i-1)>n$
vertices in $Y_k\setminus \{f(v_1),\ldots,f(v_{i-1})\}$ which are
adjacent to $f(v_j)$ for all neighbors $v_j$ of $v_i$ with $j<i$.
Since any of these vertices can be chosen for $f(v_i)$, this
completes the proof by induction.
\end{proof}

\section{Random graphs}\label{randomnext}

In this section we discuss arrangeability of sparse random graphs.
Our results imply linear upper bounds on Ramsey numbers of these
graphs. We start the section with two simple lemmas relating
$(d,\Delta)$-degeneracy with $p$-arrangeability.

\begin{lemma}\label{arrangedeg}
If a graph is $(d,\Delta)$-degenerate, then it is $\left(\Delta(d-1)+1\right)$-arrangeable.
\end{lemma}
\begin{proof}
Let $v_1,\ldots,v_n$ be a $(d,\Delta)$-degenerate ordering of the
vertices of a graph $G$. Then for each $i$, there are at most
$\Delta$ subsets $S$, each of cardinality at most $d$, such that
$S=N(v_j) \cap \{v_1,\ldots,v_i\}$ for some neighbor $v_j$ of $v_i$
with $j>i$. Therefore, for any vertex $v_i$, its neighbors to the
right of $v_i$ have together at most $\Delta(d-1)+1$ neighbors to
the left of $v_i$ (including $v_i$), and so the graph $G$ is
$\left(\Delta(d-1)+1\right)$-arrangeable.
\end{proof}

\begin{lemma}\label{arrangedeg1}
If a graph is $p$-arrangeable, then it is $(p,2^{p-1})$-degenerate.
\end{lemma}
\begin{proof}
Let $v_1,\ldots,v_n$ be a $p$-arrangeable ordering of the vertices
of a graph $G$. For any vertex $v_i$, its neighbors to the right of
$v_i$ have together at most $p$ neighbors to the left of $v_i$
(including $v_i$). Let $N(j,i)$ denote the set $v_k$ of neighbors of
$v_j$ with $k \leq i$. For every neighbor $v_j$ of $v_i$ with $j>i$,
the set $N(j,i)$ lies in a set of size $p$ that contains $v_i$, so
there are at most $2^{p-1}$ such sets $N(j,i)$. Let $v_i$ be the
neighbor of $v_j$ which has maximum index $i<j$. Then using the
$p$-arrangeability property for $v_i$, we get that the number of
neighbors of $v_j$ in $\{v_1,\ldots,v_j\}$ is $|N(j,i)| \leq p$.
Hence, the ordering demonstrates that $G$ is
$(p,2^{p-1})$-arrangeable.
\end{proof}

We prove for $d \geq 10$ that a.s.\ the random graph $G(n,d/n)$ is
$256d^2$-arrangeable. This follows from Lemma \ref{arrangedeg} and
Theorem \ref{cool} below which says that for $d \geq 10$ a.s.\
$G(n,d/n)$ is $(16d,16d)$-degenerate. Theorems \ref{cool} and
\ref{GRRtheorem} together imply Theorem \ref{corrandom1}, which says
that $G(n,d/n)$ almost surely has Ramsey number at most $2^{cd\log^2
d}n$, where $c$ is an absolute constant. Since the random bipartite
graph $G(n,n,d/n)$ is a subgraph of $G(2n,d/n)$, Theorem \ref{cool}
implies that almost surely $G(n,n,d/n)$ is $(32d,32d)$-degenerate.
Together with Theorem \ref{ramseydegenerategeneral}, this proves
Theorem \ref{corrandom2}, which says that $G(n,n,d/n)$ almost surely
has Ramsey number at most $2^{cd}n$, where $c$ is an absolute
constant.

The ordering of the vertices of $G(n,d/n)$ used to prove Theorem
\ref{cool} is a careful modification of the ordering by decreasing
degrees. Let $A$ be the set of vertices of degree more than $16d$.
It is easy to show that a.s.\ $A$ is quite small. We then enlarge
$A$ to a set $F(A)$ that we show has the property that no vertex in
the complement of $F(A)$ has more than one neighbor in $F(A)$ and
a.s.\ $|F(A)| \leq 4|A|$. Since $|F(A)|$ is small enough, a.s.\ any
set with size $|F(A)|$ (so, in particular, the set $F(A)$ itself) is
sparse enough that the subgraph of $G(n,d/n)$ induced by it is
$(2,3)$-degenerate. We first order the set $F(A)$ and then add the
remaining vertices of $G(n,d/n)$ arbitrarily. We use this vertex
order to demonstrate that a.s.\ $G(n,d/n)$ is
$(16d,16d)$-degenerate.

Before proving Theorem \ref{cool}, we need several simple lemmas.

\begin{lemma}\label{first1}
Almost surely there are at most $2^{4-8d}n$ vertices of $G(n,d/n)$
with degree larger than $16d$.
\end{lemma}
\begin{proof}
Let $A$ be the subset of $s=2^{4-8d}n$ vertices of largest degree in
$G=G(n,d/n)$ and $D$ be the minimum degree of vertices in $A$. So
there are at least $sD/2$ edges that have at least one vertex in
$A$. Consider a random subset $A'$ of $A$ with size $|A|/2$. Every
edge which contains a vertex of $A$ has a probability at least $1/2$
of having exactly one vertex in $A'$. This can be easily seen by
considering the cases that the edge lies entirely in $A$ and that
the edge has exactly one vertex in $A$. So there is a subset $A'
\subset A$ of size $|A|/2$ such that the number $m$ of edges between
$A'$ and $V(G) \setminus A'$ satisfies $m \geq sD/4=|A'|D/2$.

We now give an upper bound on the probability that $D \geq 16d$.
Each subset $A'$ of $G(n,d/n)$ of size $s/2$ has probability at most
$${\frac{s}{2}(n-\frac{s}{2}) \choose m}(d/n)^{m} \leq
\left(\frac{esn}{2m}\right)^m(d/n)^m \leq
\left(\frac{2sd}{m}\right)^m \leq \left(\frac{8d}{D}\right)^m \leq
2^{-4ds}$$ of having at least $m \geq (s/2)(16d)/2=4sd$ edges
between $A'$ and $V(G) \setminus A'$. Therefore the probability that
there is a subset $A'$ of size $s/2$ which has at least $4sd$ edges
between $A'$ and $V(G) \setminus A'$ is at most
$${n \choose s/2}2^{-4ds}
<\left(\frac{2en}{s}\right)^{s/2}2^{-4ds} \leq \left(\frac{2^{3-8d}n}{s}\right)^{s/2}=o(1),$$
completing the proof.
\end{proof}

\begin{lemma}\label{third3}
If a graph $G=(V,E)$ of order $n$ has less than $\frac{9}{8}n$ edges, then it contains a vertex of degree at most one or contains a vertex
of degree two whose both neighbors have degree two.
\end{lemma}
\begin{proof}
Suppose for contradiction that $G$ has minimum degree at least $2$
and there is no vertex of degree $2$ whose both neighbors have
degree two. Let $V_1 \subset V$ be those vertices with degree $2$
and $V_2=V \setminus V_1$ be those vertices of degree at least $3$.
Let $x=|V_1|$. Since every vertex in $V_1$ has degree at least two
and every vertex in $V_2$ has degree at least three, then the number
of edges of $G$ is at least
$\frac{2x+3(n-x)}{2}=\frac{3n}{2}-\frac{x}{2}$. Since we assumed
that every vertex of degree two has at most one neighbor with degree
two, then the subgraph of $G$ induced by $V_1$ has maximum degree at
most one. Therefore, $V_1$ spans at most $x/2$ edges. Since the
vertices in $V_1$ have degree $2$, then the number of edges of $G$
with at least one vertex in $V_1$ is $2x-e(V_1) \geq \frac{3}{2}x$.
Hence, the number of edges of $G$ is at least
$\max\left(\frac{3n}{2}-\frac{x}{2},\frac{3}{2}x \right)$.
Regardless of the value of $x$, this number is always at least
$\frac{9n}{8}$, a contradiction.
\end{proof}

\begin{lemma}\label{fourth4}
If for $r,s \geq 1$, every subgraph of a graph $G=(V,E)$ has a
vertex with degree at most one or a vertex with degree at most $s$
all of whose neighbors have degree at most $r$, then $G$ is
$(s,r+1)$-degenerate.
\end{lemma}
\begin{proof}
Pick out vertices $v_n,v_{n-1},\ldots,v_1$ one by one such that
 for each $j$, in the subgraph of $G$ induced by
$V \setminus \{v_{n},v_{n-1},\ldots,v_{j+1}\}$, the vertex $v_j$ has
degree at most one or it has degree at most $s$ and all of its
neighbors have degree at most $r$. Let $L_j=\{v_1,\ldots,v_j\}$.
Note first that this ordering has the property that each vertex
$v_j$ has at most $s$ neighbors $v_i$ with $i<j$ since its degree in
the subgraph of $G$ induced by $L_j$ is at most $s$. Let $N_1(v_j)$
be those vertices $v_k$ with $k>j$ that are adjacent to $v_j$ and
have a neighbor in $L_{j-1}$. The cardinality of $N_1(v_j)$ is at
most $r$ since otherwise the vertex $v_h \in N_1(v_j)$ with the
largest index $h$ has at least two neighbors in $L_{h}$ and has a
neighbor $v_j \in L_h$ which has degree more than $r$ in $L_{h}$,
contradicting how we chose $v_h$. The vertices $v_k$ with $k>j$ that
are adjacent to $v_j$ and have degree one in the subgraph of $G$
induced by $L_k$ satisfy $N(v_k) \cap L_j = \{v_j\}$. Therefore, for
each $j$, there are at most $r+1$ sets $S \subset L_j$ for which
 $S=N(v_k) \cap L_j$ for some vertex $v_k$ adjacent
to $v_j$ with $k>j$. Hence, this ordering shows that $G$ is
$(s,r+1)$-degenerate.
\end{proof}

\begin{lemma}\label{fifth5}
Almost surely every subgraph $G'$ of $G(n,d/n)$ with $t\leq
(5d)^{-9}n$ vertices has average degree less than $9/4$.
\end{lemma}
\begin{proof}
Let $S$ be a subset of size $t$ with $t \leq (5d)^{-9}n$. The
probability that $S$ has at least $m=\frac{9}{8}t$ edges is at most
${{t \choose 2} \choose m}(d/n)^{m}$. Therefore, by the union bound,
the probability that there is a subset of size $t$ with at least
$m=\frac{9}{8}t$ edges is at most
\begin{eqnarray*}
{n \choose t}{{t \choose 2} \choose m}\left(\frac{d}{n}\right)^m
&\leq& \left(\frac{en}{t}\right)^t \left(\frac{et^2}{2m}\right)^{m} \left(\frac{d}{n}\right)^m =
e^t \left(\frac{4e}{9}\right)^{9t/8} \left(\frac{n}{t}\right)^t \left(\frac{dt}{n}\right)^{9t/8}\\
&\leq& 5^t \left(\frac{n}{t}\right)^t \left(\frac{dt}{n}\right)^{9t/8}
=5^t(d^9t/n)^{t/8}.
\end{eqnarray*}
Summing over all $t \leq
(5d)^{-9}n$, one easily checks that the probability that there is an
induced subgraph with at most $(5d)^{-9}n$  vertices and average
degree at least $9/4$ is $o(1)$, completing the proof.
\end{proof}

For a graph $G=(V,E)$ and vertex subset $S \subset V$, let $F(S)$
denote a vertex subset formed by adding vertices from $V \setminus
S$ with at least two neighbors in $S$ one by one until no vertex in
$V \setminus F(S)$ has at least two neighbors in $F(S)$. It is not
difficult to see that $F(S)$ is uniquely determined by $S$.

\begin{lemma}\label{second2}
Almost surely every vertex subset $S$ of $G(n,d/n)$ with cardinality
$t \leq (5d)^{-10}n$ has $|F(S)|\leq 4t$.
\end{lemma}
\begin{proof}
Suppose not, and consider the set $T$ of the first $4t$ vertices of
$F(S)$. By definition, the number of edges in the induced subgraph
by $T$ is at least $2(|T|-|S|) \geq \frac{3}{2}|T|$, so the average
degree of this induced subgraph is at least $3$. But Lemma
\ref{fifth5} implies that a.s. every induced subgraph of $G(n,d/n)$
with at most $4t$ vertices has average degree less than $3$, a
contradiction.
\end{proof}

\begin{theorem}\label{cool} For $d \geq 10$, the graph $G(n,d/n)$ is almost surely $(16d,16d)$-degenerate.
\end{theorem}
\begin{proof}
Let $A$ be the $(5d)^{-10}n$ vertices of largest degree in
$G(n,d/n)$. Since $(5d)^{-10}n \geq 2^{4-8d}n$ for $d \geq 10$, by
Lemma \ref{first1}, a.s.\ all vertices not in $A$ have degree at
most $16d$. By Lemma \ref{second2}, a.s.\ $|F(A)| \leq 4|A|$. By
Lemmas \ref{fifth5}, \ref{third3}, and \ref{fourth4}, a.s. the
subgraph of $G$ induced by $F(A)$ is $(2,3)$-degenerate. Let
$v_1,\ldots,v_{|F(A)|}$ be an ordering of $F(A)$ that respects the
$(2,3)$-degeneracy. Arbitrarily order the vertices not in $F(A)$ as
$v_{|F(A)|+1},\ldots,v_n$. Let $L_{j}=\{v_1,\ldots,v_j\}$. We claim
that this is the desired ordering for the vertices of $G(n,d/n)$.
Consider a vertex $v_i$. If $i \leq |F(A)|$, then $v_i$ has the
following three properties:
\begin{itemize}
\item there are at most two neighbors $v_j$ of $v_i$ with $j<i$,
\item there are at most three subsets $S \subset L_i$ for which
there is a neighbor $v_h$ of $v_i$ with  $i<h \leq |F(A)|$ and
$N(v_h) \cap L_i=S$, and
\item every neighbor $v_k$ of $v_i$ with $k>|F(A)|$ has $N(v_k) \cap
L_i=\{v_i\}$.
\end{itemize}
If $i > |F(A)|$, then $v_i$ has maximum degree at most $16d$.
Therefore, this ordering demonstrates that $G(n,d/n)$ is a.s.\
$(16d,16d)$-degenerate, completing the proof.
\end{proof}

We next prove Lemma \ref{eight8}, which completes the proof of
Theorem \ref{randomarrange}, and says that  for $d \geq 300$ a.s.\
$G(n,d/n)$ is not $d^2/144$-arrangeable.

\begin{lemma}\label{six6}
Let $p=d/n$ with $d \geq 300$. Almost surely every pair $A,B$ of
disjoint subsets of $G(n,p)$ of size at least $n/6$ have at least
$p|A||B|/2$ edges between them.
\end{lemma}
\begin{proof}
Let $s=n/6$. If there are disjoint subsets $A,B$ each of size at
least $s$ which have less than $p|A||B|/2$ between them, then by
averaging over all subsets $A' \subset A$ and $B' \subset B$ with
$|A'|=|B'|=s$, there are subsets $A' \subset A$ and $B' \subset B$
with $|A'|=|B'|=s$ and the number of edges between $A'$ and $B'$ is
at most $p|A'||B'|/2$.

Using the standard Chernoff bound (see page 306 of \cite{AlSp}) for
$s^2$ independent coin flips each coming up heads with probability
$p=d/n$, the probability that a fixed pair of disjoint subsets
$A',B'$ of $G(n,d/n)$ each of size $s$ have less than $ps^2/2$ edges
between them is at most
$$e^{-(ps^2/2)^2/(2ps^2)}=e^{-ps^2/8}.$$

The probability that no pair of disjoint subsets $A',B'$ each of
size $s$ have less than $ps^2/2$ edges between them is at most
$${n \choose s}{n-s \choose s}e^{-ps^2/8} \leq
\left(\frac{en}{s}\right)^{2s} e^{-ps^2/8}=\left(6e \cdot e^{-d/96}
\right)^{2s}=o(1).$$ So a.s.\ every pair of disjoint subsets $A,B$
each with cardinality at least $n/6$ have at least $p|A||B|/2$ edges
between them.
\end{proof}

\begin{lemma}\label{seven7}
In $G(n,d/n)$ with $d=o(n^{1/6})$, a.s.\ no pair of vertices have
three common neighbors.
\end{lemma}
\begin{proof}
A pair of vertices together with its three common neighbors form the
complete bipartite graph $K_{2,3}$. Let us compute the expected
number of $K_{2,3}$ in $G(n,d/n)$. We pick the vertices and then
multiply by the probability that they form $K_{2,3}$. So the
expected number of $K_{2,3}$ is
$${n \choose 2}{n-3 \choose 3}(d/n)^{6} \leq n^5(d/n)^6=o(1),$$
which implies the lemma.
\end{proof}
\begin{lemma}\label{eight8}
For $d \geq 300$, almost surely $G(n,d/n)$ is not
$d^2/144$-arrangeable.
\end{lemma}
\begin{proof}
Suppose for contradiction that $G(n,d/n)$ is $d^2/144$-arrangeable,
so there is a corresponding ordering $v_1,\ldots,v_n$ of its
vertices. Let $V_1=\{v_i\}_{i \leq n/3}$, $V_2=\{v_i\}_{n/3<i \leq
2n/3}$, and $V_3 = \{v_i\}_{2n/3<i \leq n}$. Delete from $V_3$ all
vertices that have fewer than $\frac{d}{n}|V_1|/2=d/6$ neighbors in
$V_1$, and let $V_3'$ be the remaining subset of $V_3$. Note that
a.s.\ $|V_3'| \geq n/3-n/6=n/6$ since otherwise $V_3 \setminus V_3'$
and $V_1$ each have cardinality at least $n/6$ and have fewer than
$\frac{d}{n}|V_3 \setminus V_3'||V_1|/2$ edges between them, which
would contradict Lemma \ref{six6}. Delete from $V_2$ all vertices
that have fewer than $\frac{d}{n}|V_3'|/2$ neighbors in $V_3'$, and
let $V_2'$ be the remaining subset of $V_2$. Note that $|V_2'| \geq
n/3-n/6=n/6$ since otherwise $V_2 \setminus V_2'$ and $V_3'$ each
have cardinality at least $n/6$ and have fewer than $\frac{d}{n}|V_2
\setminus V_2'||V_3'|/2$ edges between them, which would contradict
Lemma \ref{six6}. Pick any vertex $v \in V_2'$. Vertex $v$ has at
least $\frac{d}{n}|V_3'|/2 \geq d/12$ neighbors in $V_3'$. Let
$U=\{u_1,\ldots,u_r\}$ denote a set of $r=d/12$ neighbors of $v$ in
$V_3'$. Let $d(u_i)$ denote the number of neighbors of $u_i$ in
$V_1$ and $d(u_i,u_j)$ denote the number of common neighbors of
$u_i$ and $u_j$ in $V_1$. Note that for each $u_i \in U$, $d(u_i)
\geq d/6$. Also, by Lemma \ref{seven7}, $d(u_i,u_j) \leq 2$ for
distinct vertices $u_i,u_j$. By the inclusion-exclusion principle,
the number of vertices in $V_1$ adjacent to at least one vertex in
$U$ is at least
$$\sum_{1 \leq i \leq r}d(u_i)-\sum_{1 \leq i < j \leq r} d(u_i,u_j)
\geq rd/6-{r \choose 2}2>d^2/144.$$ Hence, a.s.\ $G(n,d/n)$ is not
$d^2/144$-arrangeable.
\end{proof}

\section{Concluding Remarks}

In this paper we proved that for fixed $d$, the Ramsey number of the
random graph $G(n,d/n)$ is a.s.\ linear in $n$. More precisely, we
showed that there are constants $c_1,c_2$ such that a.s.\
$$2^{c_1d}\,n \leq r\left(G(n,d/n)\right) \leq 2^{c_2d \log^2 d}\, n.$$
We think that the upper bound can be further improved and that the
Ramsey number of the random graph $G(n,d/n)$ is a.s.\ at most
$2^{cd}n$ for some constant $c$.

There are many results demonstrating that certain parameters of
random graphs are highly concentrated (see, e.g., the books
\cite{Bo,JaLuRu}). Probably the most striking example of this
phenomena is a recent result of Achlioptas and Naor \cite{AcNa}.
Extending earlier results from \cite{ShSp,Lu}, they demonstrate that
for fixed $d>0$, a.s.\ the chromatic number of the random graph
$G(n,d/n)$ is $k_d$ or $k_d+1$, where $k_d$ is the smallest integer
$k$ such that $d<2k\log k$. We don't think the Ramsey numbers of
random graphs are nearly as highly concentrated. However, it would
be interesting to determine if there is a constant $c_d>0$ for each
$d>1$ such that the random graph $G(n,d/n)$ a.s.\ has Ramsey number
$\left(c_d+o(1)\right)n$.

\end{document}